\documentclass{article}
\usepackage{amssymb,amsmath,amsthm,graphicx,kotex}

\def\qed{\hfill {\hbox{${\vcenter{\vbox{               
   \hrule height 0.4pt\hbox{\vrule width 0.4pt height 6pt
   \kern5pt\vrule width 0.4pt}\hrule height 0.4pt}}}$}}}

\def\utr{\, \underline{\ast}\, }
\def\otr{\, \overline{\ast}\, }
\def\tr{\, \triangleright\, }

\theoremstyle{definition}
\newtheorem{example}{Example}
\newtheorem{definition}{Definition}

\date{}

\title{\Large \textbf{A Quick Note On Homsets and Diagrams}}

\author{Sam Nelson\footnote{Email: Sam.Nelson@cmc.edu. Partially supported by Simons Foundation collaboration grant 702597.}}

\begin{document}
\maketitle

\begin{abstract}
The homset invariant of a knot or link $L$ with respect to an algebraic 
knot coloring structure $X$ can be identified with a set of colorings
of a diagram of $L$ by elements of $X$ via an identification of diagrammatic 
generators with algebraic generators. In some cases, particularly when
either $L$ or $X$ has a high degree of symmetry, distinct homset elements 
can be represented by superficially similar $X$-colored diagrams, potentially
leading to confusion. In this short note we examine and attempt to
clarify this phenomenon. 
\end{abstract}

\parbox{5.5in} {\textsc{Keywords:} Homset invariants, Diagrammatics, Bikei invariants

\smallskip

 }

\section{Introduction}\label{I}

In this short note we aim to clear up a potential source of confusion
regarding representation and computation of homset invariants of knots and
generalized knots (including surface-links, spatial graphs, etc.) with respect 
to finite algebraic knot coloring structures via diagrams. It is sometimes 
asserted that only the cardinality of the quandle homset (for example), not 
the homset itself, is invariant, since changing the knot diagram changes the 
way the homset is represented. 

Knot theoretic coloring structures such as keis, quandles, biquandles, bikeis,
tribrackets, etc. have axioms that are chosen in order to guarantee that
the homset is invariant. More precisely, the axioms of such structures
are chosen such that Reidemeister moves on diagrams (or other combinatorial
moves depending on the type of knotted object) induce Tietze moves on the
presentation and hence induce isomorphisms of the structure. The abstract 
homset neither knows nor cares what notation we use to represent it, and 
by its construction it is clearly invariant. 

Our choice of diagram corresponds to a presentation, and the way we represent
a homset element changes with the diagram. In particular, in some cases
the same homset element can appear to have different representations on the 
same diagram due to symmetries of the knot, the coloring object or both. 
Equivalently, the same coloring of a diagram may represent different homset 
elements depending on our choice of identification of generators with 
diagrammatic elements. In this paper we will examine this situation carefully 
and aim to clear up any potential confusion.

This paper is too short to justify dividing 
into sections; it consists of a brief warm-up example followed by a few
examples to illustrate the situation along with their resolution. 
We note also that Mai Sato has explored this same topic in much more
detail in \cite{S} around the same time that this short note was written.

This paper, including all text, diagrams and computational code,
was produced strictly by the author without the use of generative AI in any 
form.

\section{Homsets and $X$-Colored Diagrams}

In recent years many invariants of knots and knotted structures (links, braids,
tangles, virtual knots, surface-links, spatial graphs, etc.) have been defined 
using \textit{homsets}, i.e., sets of homomorphisms in an appropriate category 
(groups, quandles, biquandles, tribrackets, etc.) from the fundamental
object of the knot to a generally finite coloring object. It is sometimes
said that only the cardinality of the homset is invariant, not the homset 
itself, but this objection primarily expresses the fact that our representation
of homset elements depends on our representation of the fundamental object
of the knot in question. That is, different diagrams represent the same 
coloring homomorphism differently.

Let us say that an algebraic category $\mathcal{C}$ is a \textit{knot coloring
category} if there is a functor $\mathcal{O}$ from a category of knots or 
generalized knots (e.g. braids, tangles, surface-links, virtual knots, spatial 
graphs etc.) whose objects are diagrams and morphisms are combinatorial move 
sequences to $\mathcal{C}$ such that for each object $K$ in our knot category,
the fundamental object $\mathcal{O}(K)$ is defined by a presentation consisting
of a diagram $D$ of $K$ with generators corresponding to parts of $D$ and relations
also determined by $D$. Examples include the knot group (e.g., \cite{R}), the knot 
quandle \cite{J}, the fundamental biquandle \cite{KR}, the Niebrzydowski tribracket 
(also known as the knot-theoretic ternary quasigroup) \cite{NIE} and many more.

Generally speaking, a homset element $f\in \mathrm{Hom}(\mathcal{O}(K),X)$
in an algebraic knot coloring category
can be represented with a \textit{coloring} of a diagram of $K$, typically 
an assignment of an element of $X$ to each portion of the diagram representing
a generator for the fundamental object $\mathcal{O}(K)$ of $K$ satisfying
the appropriate relations. Examples include the venerable quandle counting 
invariant, the biquandle counting invariant and the tribracket counting invariant
as well as more recent examples. 

The relevant principle is the same as the basic theorem of linear algebra, ``A 
linear transformation is determined by its values on a basis''. With a knot's
fundamental object defined via a presentation with generators and relations,
an assignment of image to each generator determines a value for every word
built from the generators, and such an assignment determines a valid 
homomorphism if and only if the relations defining the representation are 
satisfied by the choice of assignment.

As a practical matter this means that we can compute homsets by finding 
a complete set of colorings of a diagram of the knot. The role of basis 
for the input space is played 
by the diagram itself (or more precisely, the portions of the diagram 
corresponding to generators, e.g., arcs, semiarcs, regions etc.) and the role
of the matrix expressing the homomorphism is played by the coloring of the
diagram. The role of change-of-basis matrices (also sometimes called 
``transition matrices'') is played by $X$-colored Reidemeister moves.

It is therefore tempting to suppose that two colorings of a diagram represent
the same homset element if and only if they are related by $X$-colored
Reidemeister moves. This is often true, but as we will see, there 
are cases -- particularly when the knot $K$, coloring object $X$, or both 
posses high degrees of symmetry -- where the statement is not exactly true, 
or at least somewhat misleading.

In the examples we will consider here, we will use the algebraic structure
known as \textit{bikei} or \textit{involutory biquandle}, but the principle is
similar for other algebraic knot coloring structures. Let us recall the
necessary definitions (see \cite{EN} for more):

\begin{definition}
A \textit{bikei} is a set $X$ with two binary operations 
$\utr,\otr:X\times X\to X$ satisfying the axioms
\begin{itemize}
\item[(i)] For all $x\in X$, we have $x\utr x=x\otr x$, which we will denote
by $w(x)$;
\item[(ii)] For all $x,y\in X$ we have
\[\begin{array}{rcl}
(x\utr y)\utr y & = & x, \\
(x\otr y)\otr y & = & x, \\
x\utr (y\otr x)  & = & x\utr y \quad \mathrm{and} \\
x\otr (y\utr x) & = & x\otr y,
\end{array}\]
and
\item[(iii)] For all $x,y,z,\in X$ we have
\[\begin{array}{rcl}
(x\utr y)\utr (z\utr y) & = & (x\utr z)\utr (y\otr z), \\
(x\utr y)\otr (z\utr y) & = & (x\otr z)\utr (y\otr z) \quad \mathrm{and}\\
(x\otr y)\otr (z\otr y) & = & (x\otr z)\otr (y\utr z).
\end{array}\]
\end{itemize}
\end{definition}

\begin{definition}
Let $K$ be an unoriented knot, link or surface-link represented by a diagram
$D$. Then the \textit{fundamental bikei} of $K$ is the set of equivalence
classes of bikei words in a set of generators corresponding to semiarcs
in $D$ modulo the equivalence relation generated by the bikei axioms
together with the pictured \textit{crossing relations}
\[\includegraphics{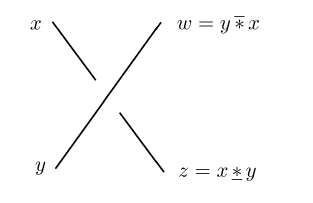}.\]
\end{definition}

Now, let us see an example to illustrate the computation of the homset as 
a set of $X$-colored diagrams.

\begin{example}
The Hopf link has fundamental bikei given by the presentation
\[\raisebox{-0.6in}{\scalebox{0.8}{\includegraphics{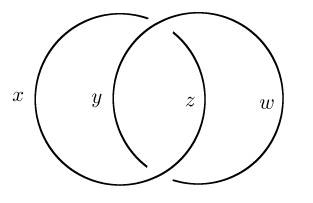}}} \quad 
\mathcal{BK}(L)=\langle x,y,z,w\ | 
z=x\utr y, z=x\otr y, w=y\utr x, w=y\otr x\rangle.\]

There are five colorings of $L$ by the bikei defined by the operation tables
\[
\begin{array}{r|rrr}
\utr & 1 & 2 & 3 \\ \hline
1 & 1 & 1 & 2 \\
2 & 2 & 2 & 1 \\
3 & 3 & 3 & 3
\end{array}
\quad \begin{array}{r|rrr}
\otr & 1 & 2 & 3 \\ \hline
1 & 1 & 1 & 1 \\
2 & 2 & 2 & 2 \\
3 & 3 & 3 & 3
\end{array}
\]
as shown
\[\includegraphics{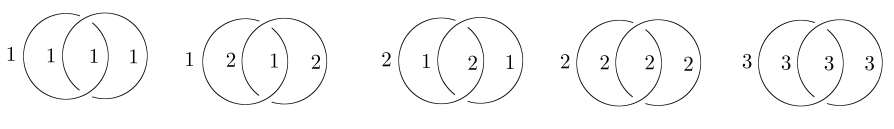}.\]
We can also think of these non-diagrammatically via a table:
\[\scalebox{0.9}{$\begin{array}{cccc|cccc}
f(x) & f(y) & f(z) & f(w) & f(z)=f(x)\utr f(y)? 
& f(z)=f(x)\otr f(y)? &  f(w)=f(y)\utr f(x)? &  f(w)=f(y)\otr f(x)? \\ \hline
1 & 1 & 1 & 1 & \checkmark & \checkmark & \checkmark & \checkmark \\
1 & 2 & 1 & 2 & \checkmark & \checkmark & \checkmark & \checkmark \\
2 & 1 & 2 & 1 & \checkmark & \checkmark & \checkmark & \checkmark \\
2 & 2 & 2 & 2 & \checkmark & \checkmark & \checkmark & \checkmark \\
3 & 3 & 3 & 3 & \checkmark & \checkmark & \checkmark & \checkmark 
\end{array}$}.
\]

Changing the diagram by a Reidemeister move changes the presentation
of the fundamental bikei by Tietze moves: 
\[\raisebox{-0.6in}{\scalebox{0.8}{\includegraphics{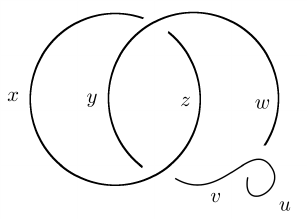}}}\quad
\begin{array}{rcl}\mathcal{BK}(L) =\langle x,y,z,w,u,v\ & | & 
z=x\utr y, z=x\otr y, v=y\utr x,  \\ 
& &  w=y\otr x, u=w\utr v, u=v\otr w\rangle.\end{array}\]
We can track the change on the $X$-colored diagrams 
\[\includegraphics{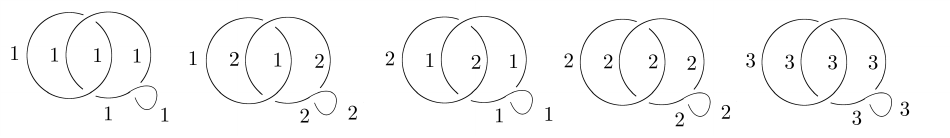}\]
or in our table:
\[\scalebox{0.85}{$\begin{array}{cccccc|ccccccc}
f(x) & f(y) & f(z) & f(w) & f(u) & f(v) & f(z)=
& f(z)= & f(v)= & f(w)=
& f(u)= & f(u)= \\ 
 &  &  &  &  &  & f(x)\utr f(y)? 
& f(x)\otr f(y)? & f(y)\utr f(x)? & f(y)\otr f(x)? 
& f(w)\utr f(v)? & f(v)\otr f(w)? \\ \hline
1 & 1 & 1 & 1 & 1 & 1 & \checkmark & \checkmark & \checkmark & \checkmark & \checkmark & \checkmark \\
1 & 2 & 1 & 2 & 2 & 2 & \checkmark & \checkmark & \checkmark & \checkmark & \checkmark & \checkmark \\
2 & 1 & 2 & 1 & 1 & 1 & \checkmark & \checkmark & \checkmark & \checkmark & \checkmark & \checkmark \\
2 & 2 & 2 & 2 & 2 & 2 & \checkmark & \checkmark & \checkmark & \checkmark & \checkmark & \checkmark \\
3 & 3 & 3 & 3 & 3 & 3 & \checkmark & \checkmark & \checkmark & \checkmark & \checkmark & \checkmark 
\end{array}$}.
\]
\end{example}

Thus, generally speaking, we can compute the homset by finding all colorings
of our choice of diagram, and changing the diagram by Reidemeister moves
yields new diagrammatic representations of the same homset elements.
What, then, can we make of the following example?

\begin{example} \label{ex2}
Consider the unknot and the bikei $X$ with operation table 
\[
\begin{array}{r|rr}
\utr & 1 & 2 \\ \hline
1 & 2 & 2 \\
2 & 1 & 1
\end{array}
\quad
\begin{array}{r|rr}
\otr & 1 & 2 \\ \hline
1 & 2 & 2 \\
2 & 1 & 1
\end{array}.
\]
This bikei can be considered as the ring $\mathbb{Z}_2$ with bikei
operations $x\utr y=x\otr y=x+1$ (where we write the class of 0 as 2). 
There are two $X$-colorings:
\[\includegraphics{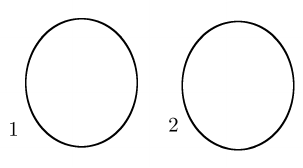}\]
But now, taking the first one, we can perform two Reidemeister I moves
\[\includegraphics{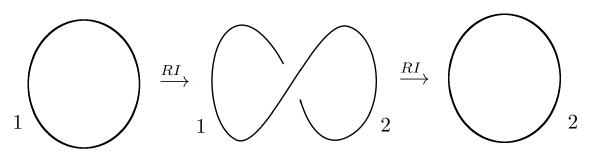}\]
and it seems we have transformed one coloring into the other! Perhaps 
there is really only one coloring? 

In fact, writing the presentations and tracking the Tietze moves, we can 
see that diagram obtained after the two RI moves still represents 
a different homomorphism from the apparently identical diagram representing
the other homset element; it still
sends the original generator $x$ to $1$ and $w(x)$ to 2.
\[\includegraphics{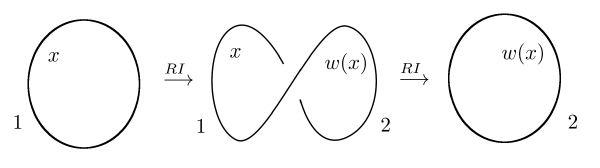}\]
\[
\begin{array}{c}
f(x) \\ \hline
1 
\end{array}
\quad\stackrel{RI}{\longrightarrow}\quad
\begin{array}{cc}
f(x) & f(w(x)) \\ \hline
1 & 2 \\
\end{array}
\quad \stackrel{RI}{\longrightarrow}\quad
\begin{array}{c}
f(w(x)) \\ \hline
 2 \\
\end{array}
\]
We can also see this by taking care to consistently apply the same
Reidemeister moves to all colorings in the homset.
\[\includegraphics{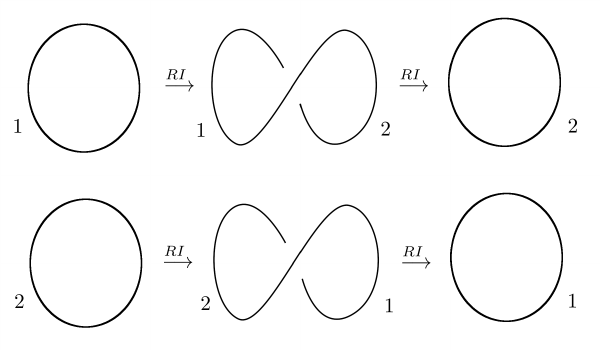}\]
\end{example}

Careful consideration of Example \ref{ex2} shows that the two apparently
identical diagrams are subtly different in ways that aren't apparent from the
diagram. For example, if we orient the knot, the two colorings by the element
$2$ have opposite orientations:
\[\includegraphics{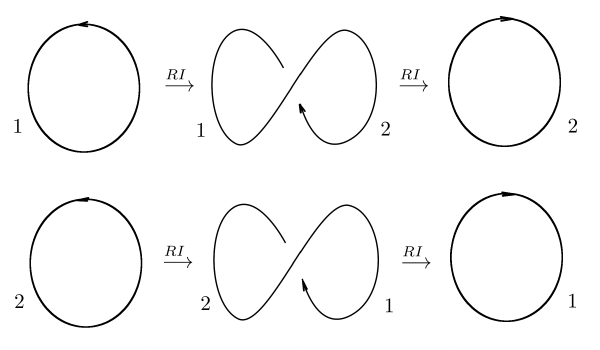}\]
Similarly, if we include a blackboard framing curve, the framing curve
flips from outside the unknot to inside:
\[\includegraphics{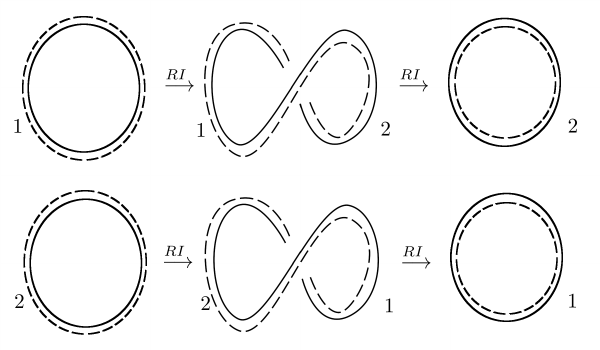}\]

\begin{example}
A similar situation happens with the colorings of the trefoil by the
bikei $X$ from Example \ref{ex2}.
\[\scalebox{1.5}{\includegraphics{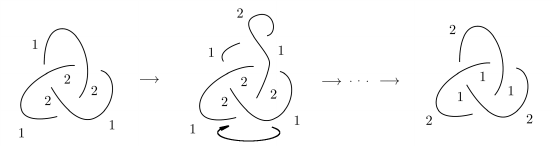}}\]
As with the unknot, this apparent ability to transform one coloring into the 
other is an illusion due to the symmetry of the diagram we are using to 
represent the trefoil.
\end{example}

\begin{example}
Even the venerable Fox colorings of the trefoil (i.e., bikei colorings by
$X=\mathbb{Z}_3$ with $x\utr y=2x+2y$ and $x\otr y=x$) have a similar kind 
of property where the six nontrivial colorings are related by a $D_3$ 
symmetry group. 
\[\scalebox{0.9}{\includegraphics{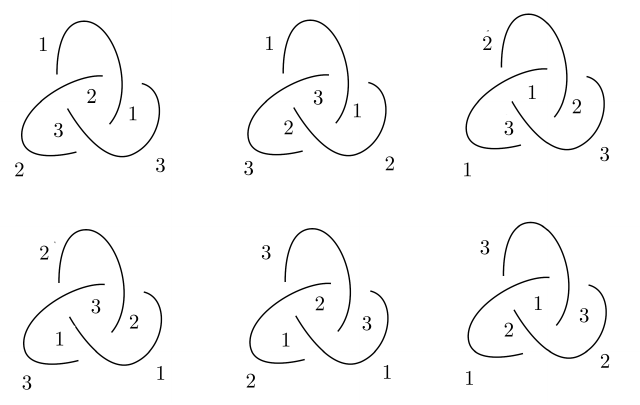}}\]
As before, the ability to match one diagram with another 
by rotation or reflection does not imply that the homomorphisms represented 
by these $X$-colored diagrams are the same; they are in this case 
sending the same generators to different elements of $X$.
\end{example}

The next example may help to clarify the situation.

\begin{example}
The unlink of two components has fundamental bikei given by the free bikei
on two generators, $\mathcal{BK}(U_2)=\langle g_1, g_2\ |\ \rangle$. A homset 
element is then represented by a choice of image for each of two circles. 
\[\scalebox{1.25}{\includegraphics{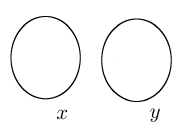}}\]
However, from a coloring alone we cannot tell which circle
represents which generator of the free bikei; is this diagram representing 
$\{f(g_1)=x, f(g_2)=y\}$ or $\{f(g_1)=y, f(g_2)=x\}$?

To clear up the potential confusion, it suffices to identify clearly 
on the diagram which diagrammatic generators (arcs, semiarcs, regions, 
etc.) are identified with which algebraic generators in our presentation.
\[\scalebox{1.25}{\includegraphics{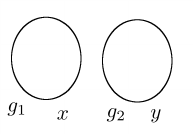}}
\quad\raisebox{0.5in}{$\ne$} \quad
\scalebox{1.25}{\includegraphics{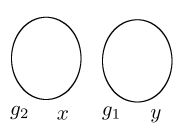}}
\]
\end{example}

Moreover, even though $X$-colorings of knot or link
diagrams are a convenient and approachable way to do homset computations,
not every presentation of a knot's fundamental bikei (or knot quandle, etc.)
is directly representable as an $X$-colored diagram. 

\begin{example}
The fundamental bikei of the trefoil knot has presentation 
\[\langle x,y\ |\ (x\utr y)\otr y=(y\otr x)\utr x,\ 
x\otr((x\utr y)\otr y)=y\utr((x\utr y)\otr y)
\rangle\]
which clearly cannot correspond to a diagram as the trefoil has crossing 
number 3 and hence no diagrams with only two semiarcs.
\end{example}

Nevertheless, by definition the fundamental object of a knot or link in an
algebraic knot coloring category has a presentation representable with
a diagram. Of course, not all objects in these categories are necessarily the
fundamental object of a knot or link. It is shown in \cite{NAGT} that every 
presentation of a finitely generated quandle, for example, can be modified 
via Tietze moves to \textit{short form}, i.e., such 
that every relation has the form $x=y\tr z$; each such relation can then be 
interpreted as a choice of one of two oriented crossings which can then be 
connected by matching colored endpoints to form a tangle such that the
fundamental quandle of the tangle is isomorphic to the presented quandle. 
Even when such a tangle has no open endpoints and thus forms a link, 
however, the resulting diagram is defined at most up to welded isotopy and
mirror image. 

We conclude, therefore, with the following observations:
\begin{itemize}
\item[(i)] In some cases, in the absence of explicit identification 
of diagrammatic generators with algebraic generators, different homset 
elements may appear to be represented by the same $X$-colored diagram,
\item[(ii)] This does not cause a problem for computation and representation
of the homset since we can work with a fixed choice of diagram, resulting
in a distinct representation for each homset element as a coloring of the
fixed choice of diagram explicitly corresponding to a presentation, and
\item[(iii)] If we need to change the diagram by combinatorial moves, we can 
track the correspondence with an algebraic presentation by simply labeling 
the diagrammatic generators with their corresponding algebraic generators.
\end{itemize}

\bibliography{sn-solo25}{}
\bibliographystyle{abbrv}

\bigskip

\noindent
\textsc{Department of Mathematical Sciences \\
Claremont McKenna College \\
850 Columbia Ave. \\
Claremont,  CA 91711}

\end{document}